\documentclass[12pt]{amsart}
\usepackage{amssymb}
\newtheorem{theorem}{Theorem}
\newtheorem{proposition}[theorem]{Proposition}
\newtheorem{corollary}[theorem]{Corollary}

\def\ds{\displaystyle}
\def\dist{\operatorname{dist}}
\def\diam{\operatorname{diam}}

\def\DD{\Bbb D}
\def\CC{\Bbb C}

\title{Lipschitzness of the Lempert and Green functions}

\author{Nikolai Nikolov, Peter Pflug and Pascal J. Thomas}

\address
{Institute of Mathematics and Informatics\\ Bulgarian Academy of
Sciences\\ Acad. G. Bonchev 8, 1113 Sofia,
Bulgaria}\email{nik@math.bas.bg}

\address{Carl von Ossietzky Universit\"at Oldenburg\\
Institut f\"ur Mathe\-ma\-tik\\ Postfach 2503\\ D-26111 Oldenburg,
Germany}\email{pflug@mathematik.uni-oldenburg.de}

\address{Laboratoire Emile Picard, UMR CNRS 5580\\
Universit\'e Paul Sabatier, 118 Route de Narbonne\\ F-31062
Toulouse Cedex, France} \email{pthomas@cict.fr}

\subjclass[2000]{32F45, 32U45.}

\keywords{Lipschitzness, Lempert function, Kobayashi--Royden
pseudometric, pluricomplex Green function, Azukawa pseudometric.}

\begin{document}

\begin{thanks}{This paper was started during the stay of the first named
author at the Carl von Ossietzky Universit\"at, Oldenburg (October
2007; supported by a grant from the DFG, Az.~PF 227/9-1) and it
was finished during his stay at the Universit\'e Paul Sabatier,
Toulouse (January 2008).}
\end{thanks}

\begin{abstract} Necessary and sufficient conditions for
Lipschitzness of the Lempert and Green functions are found in
terms of their boundary behaviors.
\end{abstract}

\maketitle

\section{Introduction and results}

By $\Bbb D$ we denote  the unit disc in $\Bbb C.$ Let $D$ be a
domain in $\Bbb C^n.$ Recall first the definitions of the Lempert
function and the Kobayashi--Royden pseudometric of $D:$
$$l_D(z,w):=\inf\{\alpha\in[0,1):\exists\varphi\in\mathcal O(\Bbb
D,D):\varphi(0)=z,\varphi(\alpha)=w\},$$
$$\kappa_D(z;X):=\inf\{\alpha\ge 0:\exists\varphi\in\mathcal O(\Bbb
D,D):\varphi(0)=z,\alpha\varphi'(0)=X\}.$$ We point out that both
functions are upper semicontinuous and $l_D$ is symmetric. The
Kobayashi--Buseman pseudometric $\hat\kappa_D(z;\cdot)$ (the
Kobayashi pseudodistance $k_D$) is the largest pseudonorm
(pseudodistance) which does not exceed $\kappa_D(z;\cdot)$
($\tanh^{-1}l_D$). Note that if $D$ is a \textit{taut} domain,
i.e., $\mathcal O(\Bbb D,D)$ is a normal family, then $\kappa_D$
and $\hat\kappa_D$ are the infinitesimal forms of $l_D$ and $k_D,$
respectively (see \cite{Nik-Pfl2}, Theorem 1 for a more general
result). Moreover, recall (cf. \cite{Nik-Pfl0}, Proposition 3.2)
that $D$ is a taut domain  if and only if
$$\lim_{z\in K,w\to\partial D}l_D(z,w)=1\mbox{ for any }K\Subset D.$$
(Note that for a unbounded $D$ the point $\infty$ belongs, by
definition, to $\partial D$.)

The main result in \cite{Kra} (see Theorem 6 there) is that
$\kappa_D$ is a locally H\"older function of order $2/3$ on any
$C^6$-smooth strongly pseudoconvex domain $D$ in $\Bbb C^n$ (see
also \cite{Kra1}, where it is claimed that $\kappa_D$ is locally
Lipschitz but the proof there seems to be non correct).

Our first goal in the present note is to generalize this result
showing that $l_D$ and $\kappa_D$ are Lipschitz functions under a
natural assumption about the boundary behavior of $l_D.$ In fact
we have the following result.\footnote{Proofs for this and the
next results will be presented in section 2.}

\begin{proposition}\label{pr1} Let $D\subset\Bbb C^n $ be a hyperbolic
domain (i.e., $k_D$ is a distance) and $K\Subset D$ be such that
$\ds\sup_{z\in K,w\in D}\frac{1-l_D(z,w)}{\dist(w,\partial
D)}<\infty.$ Then:

{\rm(i)} $l_D$ is a Lipschitz function on $K\times D;$

{\rm(ii)} there is a $C>0$ such that if $z,w\in K,$ $X,Y\in\Bbb
C^n,$ then
$$|\kappa_D(z;X)-\kappa_D(w;Y)|\le C((\|X\|+\|Y\|)\cdot\|z-w\|+\|X-Y\|).$$
\end{proposition}

\noindent{\bf Remark A.} (a) By symmetry, $l_D$ is a Lipschitz
function on $D\times K,$ too. On the other hand, $l_{\Bbb D}$ is
not a Lipschitz function on $\DD\times \DD.$

(b) For the Carath\'eodory--Reiffen pseudometric the same estimate
as in (ii) remains true for any domain in $\Bbb C^n$ (see
\cite{Jar-Pfl1}, Proposition 2.5.1(c)).

(c) Any compact subset of a strongly pseudoconvex domain satisfies
the assumption of Proposition \ref{pr1} (cf. \cite{Jar-Pfl1},
Theorem 10.2.1).

(d) If $D\subset\Bbb C^n $ is a hyperbolic domain, $K\Subset D,$
$L\Subset\Bbb C^n$, and $\ds\sup_{z\in K,w\in
D}\frac{1-l_D(z,w)}{\dist^\alpha(w,\partial D)}<\infty,$
$\alpha\in(0,1),$ then obvious modifications in the proof of
Proposition \ref{pr1} imply that $l_D$ and $\kappa_D$ are H\"older
functions with exponent $\alpha$ on $K\times D$ and $K\times L,$
respectively. On the other hand, $\alpha$ cannot be taken larger
than 1; one can show that for any domain $D\subsetneq\Bbb C^n$ and
any point $z\in D$ we have $\ds\limsup_{w\to\partial
D}\frac{1-l_D(z,w)}{\dist(w,\partial D)}>0.$
\smallskip

We point out that for a taut domain $D$ the assumption of
Proposition \ref{pr1} is also necessary for $l_D$ to be a
Lipschitz function.

\begin{corollary}\label{cor2} Let $D\subset\Bbb C^n $ be a taut domain
and $K\Subset D.$ Then $\ds\sup_{z\in K,w\in
D}\frac{1-l_D(z,w)}{\dist(w,\partial D)}<\infty$ if and only if
$l_D$ is a Lipschitz function on $K\times D.$
\end{corollary}

To prove the Lipschitzness of $\hat\kappa_D$ under the assumption
of Proposition \ref{pr1}, we shall need the following result.

\begin{proposition}\label{pr3} Let $D\subset\Bbb C^n$ be a
hyperbolic domain and let $K\Subset D,$ $c>0$ be such that
$$|\kappa_D(z;X)-\kappa_D(w;X)|\le c\|X\|\cdot \|z-w\|,\
\quad z,w\in K,\;X\in\Bbb C^n.$$ Then there is a $C>0$ such that
if $z,w\in K,$ $X,Y\in\Bbb C^n,$ then
$$|\hat\kappa_D(z;X)-\hat\kappa_D(w;Y)|\le C((\|X\|+\|Y\|)\cdot\|z-w\|+\|X-Y\|).$$
\end{proposition}

The next corollary is an immediate consequence of Propositions
\ref{pr1} and \ref{pr3}.

\begin{corollary}\label{cor4} Let $D\subset\Bbb C^n $ and $K\Subset D$ be
as in Proposition \ref{1}. Then there is a $C>0$ such that if
$z,w\in K,$ $X,Y\in\Bbb C^n,$ then
$$|\hat\kappa_D(z;X)-\hat\kappa_D(w;Y)|\le C((\|X\|+\|Y\|)\cdot\|z-w\|+\|X-Y\|).$$
\end{corollary}

The second aim of our paper is to find a necessary and sufficient
condition for the exponential of the pluricomplex Green function
to be Lipschitz (similar to that for the Lempert function).

Recall first the definitions of the pluricomplex Green function
and the Azukawa pseudometric of a domain $D$ in $\Bbb C^n:$
\begin{multline*}
g_D(z,w):=\sup\{u(w):u\in PSH(D),u<0,\\
\limsup_{\zeta\to z} (u(\zeta)-\log\|\zeta-z\|)<\infty\},
\end{multline*}
$$A_D(z;X):=\limsup_{t\nrightarrow 0}\frac{\tilde
g_D(z,z+tX)}{|t|},$$ where $\tilde g_D:=\exp g_D.$ We point out
that both functions are upper semicontinuous (cf. \cite{Jar-Pfl2},
page 10) and $\tilde g_D\le l_D.$ Note also that, in general,
$g_D$ is not symmetric.

Recall also that a domain $D\subset\Bbb C^n$ is called
\textit{hyperconvex} if it has a negative plurisubharmonic
exhaustion function. The next proposition is a consequence of the
proof of Theorem 3.1 in \cite{Blo2} (see also \cite{Blo1}, Theorem
2 for a weaker version).

\begin{proposition}\label{blo} Let $D\subset\Bbb C^n$ be a bounded
domain. Then the following conditions are equivalent:

(i) there is $u\in PSH(D)$ with $u<0$ and $\ds\inf_{z\in D
}u(z)/\dist(z,\partial D)>-\infty;$

(ii) $D$ is hyperconvex and there are $z_0\in D$ and $C>0$ such
that if and $w_1,w_2\in D\setminus\{z_0\},$ then
$$|g_D(z_0,w_1)-g_D(z_0,w_2)|\le C\frac{\|w_1-w_2\|}{\min\{\|z_0-w_1\|,\|z_0-w_2\|\}};
$$

(iii) $D$ is hyperconvex and for any $K\Subset D$ there is a $C>0$
such that if $z\in K$ and $w_1,w_2\in D\setminus\{z\},$ then
$$|g_D(z,w_1)-g_D(z,w_2)|\le C\frac{\|w_1-w_2\|}{\min\{\|z-w_1\|,\|z-w_2\|\}}.
$$
\end{proposition}

As a simple consequence we get the following result for $\tilde
g_D$.

\begin{corollary}\label{cor6} Let $D$ and $u$ be as
in Proposition \ref{blo}(i) and let $K\Subset D.$ Then there is
$C>0$ such that
$$|\tilde g_D(z,w_1)-\tilde g_D(z,w_2)|\le C\|w_1-w_2\|,\quad
z\in K,\;w_1,w_2\in D.$$
\end{corollary}

\noindent{\bf Remark B.} Let $D$ be a hyperconvex domain (not
necessary bounded) and $u$ be as in Proposition \ref{blo} (if $D$
is bounded, then (i) implies that $u$ is an exhaustion function of
$D$ and hence $D$ is hyperconvex). Then, for an arbitrary
$K\Subset D$, the assumptions of Proposition \ref{pr1} are
satisfied. Indeed, it follows from (\ref{5}) below that
\begin{multline*}
1-l_D(z,w)\leq 1-\tilde g_D(z,w)\leq-g_D(z,w)\leq c\dist(w,\partial D),\\
\quad z\in K,\; w\in D, \text{ near } \partial D;
\end{multline*}
hence the inequality in the assumption of Proposition \ref{pr1} is
fulfilled. It remains to use that $l_D\ge\tilde g_D$ and that $D$
is hyperconvex. Hence $D$ is taut (cf. \cite{Nik-Pfl0}, page 607)
and therefore hyperbolic.
\smallskip

>From Corollary \ref{cor6} we get that under the same assumptions
$\tilde g_D$ and $A_D$ are Lipschitz functions (in both
arguments).

\begin{proposition}\label{pr7} Let $D$  and $u$ be as
in Proposition \ref{blo}(i) and let $K\subset D$ be compact. Then:

(i) $\tilde g_D$ is a Lipschitz function on $K\times D;$

(ii) there is a $C>0$ such that if $z,w\in K,$ $X,Y\in\Bbb C^n,$
then
$$|A_D(z;X)-A_D(w;Y)|\le C((\|X\|+\|Y\|)\cdot\|z-w\|+\|X-Y\|).$$
\end{proposition}

It remains an open question whether $\tilde g_D$ is a Lipschitz
function on $D\times K.$
\smallskip

\noindent{\bf Remark C.} Let $D\subset\Bbb C^n$ be a pseudoconvex
balanced domain with Minkowski function $h_D.$ Recall that (cf.
\cite{Jar-Pfl1}, Propositions 3.1.10 and 4.2.7 (b))
$$l_D(0,\cdot)=\kappa_D(0;X)=g_D(0,\cdot)=A_D(0;\cdot)=h_D.$$
Note also that (cf. \cite{Nik-Pfl0}, Proposition 4.4.

\noindent D is taut $\Leftrightarrow$ D is hyperconvex
$\Leftrightarrow$ D is bounded and $h_D$ is continuous.

By Corollary \ref{cor2} or Corollary \ref{cor6}, for a taut
balanced domain $D$ the following are equivalent:

(i) there is $c>0$ such that $1-h_D(z)\le c\cdot\dist(z,\partial
D),$ $z\in D;$

(ii) there is $c'>0$ such that $|h_D(z)-h_D(w)|\le c'\|z-w\|,$
$z,w\in D.$

(Taking $h_D(z)=|z_1|+|z_2|+\sqrt{|z_1z_2|}$ provides an example
of a taut balanced domain $D\subset\Bbb C^2$ which does not have
the above properties.)

We point out that $(i)\Leftrightarrow(ii)$ with $c=c'$ for any
balanced domain $D$ in $\Bbb C^n.$

Indeed, assume that $(i)$ holds. Then for any $z,w$ with
$1>h_D(z)>h_D(w)$ we have
$$h_D(z)-h_D(w)=h_D(z)(1-h_D(w/h_D(z)))\le h_D(z)c\cdot\dist(w/h_D(z),\partial D)$$
$$\le h_D(z)c\|w/h_D(z)-z/h_D(z)\|=c\|z-w\|.$$

Conversely, assume that (ii) is true. Fix a $z\in D$. If
$\|u\|<r_z:=(1-h_D(z))/c'$, then $h_D(z+u)\le h_D(z)+c'\|u\|<1$,
which shows that $\Bbb B_n(z,r_z)\subset D.$ Hence
$\dist(z,\partial D)\ge r_z,$ that is, (i) holds with $c=c'.$

\section{Proofs}

{\it Proof of Proposition \ref{pr1}.} The assumption of
Proposition \ref{pr1} means that there is a $c>0$ such that for
any $r\in(0,1)$ and $\varphi\in\mathcal O(\Bbb D,D)$ with
$\varphi(0)\in K$ one has that
$$c\cdot\dist(\varphi(r\Bbb D),\partial D)\ge 1-r.$$

Note that there is a $c_1>0$ such that
$$c_1\|z-w\|\ge l_D(z,w),\quad z\in K,w\in D.$$
On the other hand, if $D$ is unbounded, then, by hyperbolicity,
$m^\ast=\liminf_{z\in K,w\to\infty}l_D(z,w)>0$ (use e.g.
\cite{Nik-Pfl0}, Proposition 3.1). Fix a $m\in
(0,\min\{1/2,m^\ast\})$. Then, again by hyperbolicity, we find a
$c_2>0$ such that:
$$l_D(z,w)\le m,z\in K,w\in D\Rightarrow l_D(z,w)\ge c_2\|z-w\|$$
(apply e.g. \cite{Jar-Pfl1}, Theorem 7.2.2; if $D$ is bounded, the
last inequality holds even on $K\times D$ with suitable $c_2>0$;
no other assumptions are needed in this situation). We may assume
that $c_1>1>c_2$. Set $c_3=c_1(1+c/(mc_2)).$ To prove (i), it
suffices to show that if
\begin{equation}\label{1}
|l_D(z,w_1)-l_D(z,w_2)|\le c_3\|w_1-w_2\|,\ z\in K,w_1,w_2\in D,
\end{equation}

\begin{equation}\label{1'}
|l_D(w_1,z)-l_D(w_2,z)|\le 2c_3\|w_1-w_2\|,\ w_1,w_2\in K,z\in D.
\end{equation}

To prove (\ref{1}), we may assume that $\alpha:=l_D(z,w_1)\le
l_D(z,w_2)$ and $z\neq w_1.$ Then, by hyperbolicity, $\alpha>0.$
Set $r=1-c\|w_1-w_2\|/\alpha.$ We shall consider three cases.

Case 1. $r>\max\{\alpha,m\}.$ Then for any $\alpha'\in(\alpha,r)$
there is $\varphi\in\mathcal O(\Bbb D,D)$ with $\varphi(0)=z$ and
$\varphi(\alpha')=w_1.$ Set
$\psi(\zeta)=\varphi(r\zeta)+(w_2-w_1)r\zeta/\alpha',$
$\zeta\in\Bbb D.$ Then $\psi\in\mathcal O(\Bbb D,D)$ and
$\psi(\alpha'/r)=w_2$ ($\alpha'<r$). It follows that
$l_D(z,w_2)\le\alpha/r$ and hence
$$l_D(z,w_2)-l_D(z,w_1)\le\alpha(1-r)/r=$$
$$c\|w_2-w_1\|/r\le c\|w_2-w_1\|/m\le c_3\|w_2-w_1\|.$$

Case 2. $\alpha\ge\max\{r,m\}.$ Then
$$l_D(z,w_2)-l_D(z,w_1)<1-\alpha\le 1-r=$$
$$c\|w_1-w_2\|/\alpha\le c\|w_1-w_2\|/m<c_3\|w_1-w_2\|.$$

Case 3. $m\ge\max\{r,\alpha\}.$ Then
$$\|w_1-w_2\|=(1-r)\alpha/c\ge(1-m)\alpha/c\ge (1-m)c_2\|z-w_1\|/c,$$
and, by the triangle inequality, $\|z-w_2\|\le
(1+c/((1-m)c_2))\|w_1-w_2\|.$ Since $m\le 1/2,$ it follows that
$$l_D(z,w_2)-l_D(z,w_1)<l_D(z,w_2)\le c_1\|z-w_2\|\le
c_3\|w_1-w_2\|.$$

This completes the proof of (\ref{1}).

The proof of (\ref{1'}) is similar to that of (\ref{1}) and we
sketch it. We may assume that $0<\beta:=l_D(w_1,z)\le l_D(w_2,z)$
and then set $s=1-2c\|w_1-w_2\|/\beta.$ We get as above that:

Case 1. If $\beta\ge\max\{s,m\},$ then
$l_D(w_2,z)-l_D(w_1,z)<2c\|w_1-w_2\|/\beta;$

Case 2. If $m\ge\max\{s,\beta\},$ then
$l_D(w_2,z)-l_D(w_1,z)<c_3\|w_1-w_2\|.$

Case 3. In the remaining case $s>\max\{\beta,m\},$ for any
$\beta'\in(\beta,s)$ we may find $\varphi\in\mathcal O(\Bbb D,D)$
with $\varphi(0)=w_1$ and $\varphi(\beta')=z_1.$ Set
$\psi(\zeta)=\varphi(s\zeta)+(w_2-w_1)(1-s\zeta/\beta'),$
$\zeta\in\Bbb D.$ Then $\psi\in\mathcal O(\Bbb D,D),$
$\psi(0)=w_2$ and $\psi(\beta'/s)=z.$ It follows that
$l_D(w_2,z)\le\beta/s$ and hence
$$l_D(w_2,z)-l_D(w_1,z)\le2c\|w_2-w_1\|$$ which completes the
proof of (\ref{1'}).

Next, we shall prove (ii). It is enough to show that
\begin{equation}\label{2}|\kappa_D(z;X)-\kappa_D(w;X)|\le 4cc_4\|X\|\cdot\|z-w\|,
\end{equation}
and
\begin{equation}\label{3}|\kappa_D(z;X)-\kappa_D(z;Y)|\le c_5\|X-Y\|,
\end{equation}
for any $z,w\in K,$ $X,Y\in\Bbb C^n,$ where $c_4:=\sup_{u\in
K,\|U\|=1}\kappa_D(u;U),$ $c_5:=c_4(1+2c/c_6)$ and
$c_6:=\inf_{u\in K,\|U\|=1}\kappa_D(u;U)$ ($c_6>0$ by
hyperbolicity; cf. \cite{Jar-Pfl1}, Theorem 7.2.2).

For proving (\ref{2}), observe that
$$|\kappa_D(z;X)-\kappa_D(w;X)|\le 2c_4\|X\|.$$
So (\ref{2}) is trivial if $p=1-c\|z-w\|\le 1/2.$ Otherwise, we
may assume that $\kappa_D(z;X)\le\kappa_D(w;X).$ For any
$\varphi\in\mathcal(\Bbb D,D)$ set
$\psi(\zeta)=\varphi(p\zeta)+w-z,$ $\zeta\in\Bbb D.$ Then
$\psi\in\mathcal O(\Bbb D,D)$ which shows that
$\kappa_D(w;X)\le\kappa_D(z;X)/p.$ This implies (\ref{2}) with
$2cc_4$ instead of $4cc_4.$

To get (\ref{3}), we may assume that $\gamma=\kappa_D(z;X)\le
\kappa_D(z;Y)$ and $X\neq 0.$ Then $\gamma>0.$ For
$q=1-c\|X-Y\|/\gamma$ we have two cases.

Case 1. $q>1/2.$ Let $\varphi\in\mathcal O(\Bbb D, D)$ be such
that $\varphi(0)=z$ and $\gamma'\varphi'(0)=X$ for some $\gamma'.$
Set $\psi(\zeta)=\varphi(q\zeta)+(Y-X)q\zeta/\gamma'.$ Then
$\psi\in\mathcal O(\Bbb D,\Bbb D)$ and $\gamma'\psi'(0)=qY.$ It
follows that $\kappa_D(z;Y)\le\gamma/q$ and hence
$$\kappa_D(z;Y)-\kappa_D(z;X)\le\gamma(1-q)/q=c\|X-Y\|/q\le
c_5\|X-Y\|.$$

Case 2. $q\le 1/2.$ Then $\|X-Y\|=(1-q)\gamma/c\ge c_6\|X\|/(2c)$
and, by the triangle inequality, $\|Y\|\le(1+2c/c_6)\|X-Y\|.$ It
follows that
$$\kappa_D(z;Y)-\kappa_D(z;X)<\kappa_D(z;Y)\le c_4\|Y\|\le
c_5\|X-Y\|.$$

This completes the proof of Proposition \ref{pr1}. \qed

\smallskip

{\it Proof of Corollary \ref{cor2}.}  By Proposition \ref{pr1}, it
is enough to show that if
$$|l_D(z,w_1)-l_D(z,w_2)|\le c\|w_1-w_2\|,\quad z\in K,\;w_1,w_2\in
D,$$ then $\ds\sup_{z\in K,w\in
D}\frac{1-l_D(z,w)}{\dist(w,\partial D)}<\infty.$ Suppose this is
not true. Then there are sequences $(z_j)_j\subset K$ and
$(w_j)_j\subset D$ such that $$1-l_D(z_j,w_j)\geq j
\dist(w_j,\partial D),\quad j\in\Bbb N.$$ Choose $b_j\in\partial
D$ with $\|w_j-b_j\|=\dist(w_j,\partial D)$ and sequences
$(b_{j,k})_k\subset D$ with $b_{j,k}\to b_j$ if $k\to\infty$,
$j\in\Bbb N$. Then
\begin{multline*}
j\leq \frac{1-l_D(z_j,b_{j,k})+l_D(z_j,b_{j,k})-l_D(z_j,w_j)}
{\|w_j-b_j\|}\\\leq
\frac{1-l_D(z_j,b_{j,k})}{\|w_j-b_j\|}+c\frac{\|w_j-b_{j,k}\|}{\|w_j-b_j\|}
\leq 1+2c,\end{multline*} if $k=k_j$ is sufficiently large. Recall
that $D$ is taut, therefore such $k_j$ always exist. A
contradiction.
\smallskip

{\it Proof of Proposition \ref{pr3}.} Since
$\hat\kappa_D(z;\cdot)$ is a norm, it is enough to show that for
any $K\Subset D$ there is $C'>0$ such that
$$|\hat\kappa_D(z;X)-\hat\kappa_D(w;X)|\le C'\|X\|\cdot \|z-w\|,\
\quad z,w\in K,X\in\Bbb C^n.$$ We may assume that $z\neq w,$
$X\neq 0$ and $\hat\kappa_D(z;X)\le\hat\kappa_D(w;X).$ Then there
are vectors $X_1,\dots,X_{2n-1}\in\Bbb C^n$ with sum $X$ such that
(see \cite{Nik-Pfl1}, Theorem 1)
$$\sum_{j=1}^{2n-1}\kappa_D(z;X_j)\le\hat\kappa_D(z;X)+\|X\|\cdot\|z-w\|.$$
It follows that
$$0\le\hat\kappa_D(w;X)-\hat\kappa_D(z;X)\le\sum_{j=1}^{2n-1}(\kappa_D(w;X_j)-
\kappa_D(z;X_j))+\|X\|\cdot\|z-w\|$$
$$\le \|z-w\|(\|X\|+c\sum_{j=1}^{2n-1}\|X_j\|).$$
It remains to use that
$\sum_{j=1}^{2n-1}\|X_j\|\le\frac{c_4}{c_6}\|X\|,$ where $c_4$ and
$c_6$ are as in the proof of Proposition \ref{pr1}.
\smallskip

{\it Proof of Proposition \ref{blo}.} (ii)$\Rightarrow$ (i). Put
$u=\tilde g_D(z_0,\cdot)-1$ Since $D$ is a hyperconvex domain,
then $\lim_{w\to\partial D} u(w)=0.$ Now a similar argument as in
the proof of Corollary 2 implies that $u$ has the required
property.

Since (iii)$\Rightarrow$ (ii) is trivial, it remains to prove:

(i) $\Rightarrow$ (iii). Since $u$ is an exhaustion function of
$D,$ it follows that $D$ is hyperconvex.

Fix a $K\Subset D$. We shall show that if $D$ is hyperconvex (not
necessary bounded) and $u$ is as in (i), then
\begin{equation}\label{5}
\liminf_{z\in K,w\to\partial D}g_D(z,w)/\dist(z,\partial
D)>-\infty.
\end{equation}
Indeed: Let $\tilde u$ be an exhaustion function of $D$ and $\hat
u=\max\{u,\tilde u\}.$ Then take a domain $G_1\Subset D$,
$K\Subset G_1$, and put $\varepsilon=\sup_{G_1}\hat u/2<0.$ Next
we choose a domain $G_2\Subset D$, $G_1\Subset G_2$, such that
$\inf_{\partial G_2}\hat u\ge\varepsilon.$

Fix a $z\in K$. Set $\varphi(z,\cdot)=\log(\|\cdot-z\|/\diam
G_2),$ $m=\inf_{K\times G_1}\varphi$ and
$$v_z=\left\{\begin{array}{ll}
\varphi(z,\cdot)+m&\mbox{on }G_1\\
\max\{\varphi(z,\cdot)+m,m\hat u/\varepsilon\}&\mbox{on }G_2\setminus G_1\\
m\hat u/\varepsilon&\mbox{on }D\setminus G_2\\
\end{array}.\right.$$
It is easy to check that $v_z\in PSH(D)$ for $z\in K.$ Hence
$g_D(z,\cdot)\ge v_z$ which implies (\ref{5}).

Let now $r>0$ be such that $\Bbb B(a,r)\subset D$ for any $z\in
K.$ For any $\varepsilon\in(0,r)$ we set
$$g_D^\varepsilon(z,w)=\sup\{u(w):u\in
PSH(D),u<0,u|_{\Bbb B(z,\varepsilon)}\le\log(\varepsilon/r)\}.$$
One can easily check that $g_D^\varepsilon(z,\cdot)$ is a maximal
plurisubharmonic function on $D\setminus\overline{\Bbb
B(z,\varepsilon)}$ (cf. \cite{Jar-Pfl1}, page 383 for this
notion),
\begin{equation}\label{6}
\max\{\log(\varepsilon/r),g_D(z,w)\} \le
g_D^\varepsilon(z,w)\le\log\frac{\max\{\|z-w\|,\varepsilon\}}{r}
\end{equation}
and $g_D^\varepsilon(z,\cdot)\downarrow g_D(z,\cdot)$ as
$\varepsilon\downarrow 0$ locally uniformly in $D\setminus\{z\}$
(cf. \cite{Blo2}, page 338 and Proposition 2.2). Moreover, since
$D$ is hyperconvex, $g_D^\varepsilon$ can be extended as a
continuous function on $D\times\overline D$ by setting
${g_D^\varepsilon}_{|D\times\partial D}=0.$

We shall find $c_1,c_2>0$ such that if $z\in K,$ $w_1,w_2\in
D\setminus\{z\}$, and $\varepsilon>0$ satisfy the inequality
\begin{equation}\label{7}
\max\{\varepsilon,c_1\|w_1-w_2\|\}<\min\{r/2,\|z-w_1\|,\|z-w_2\|\},
\end{equation}
then
\begin{equation}\label{8}
|g^\varepsilon_D(z,w_1)-g^\varepsilon_D(z,w_2)|\le
c_2\frac{\|w_1-w_2\|}{\min\{\|z-w_1\|,\|z-w_2\|\}}
\end{equation}

Assuming (\ref{8}), take arbitrary points $w_1,w_2\in
D\setminus\{z\}.$ To prove (iii), we may assume that
$g^\varepsilon_D(z,w_1)\le g^\varepsilon_D(z,w_2),$ where
$\varepsilon$ is as above. There is a semicircle with diameter
$[w_1w_2],$ say $\gamma:[0,\pi]\to\Bbb C^n,$
$\gamma(0)=w_1,\gamma(\pi)=w_2,$ such that
$dist(z,\gamma)=\min\{\|z-w_1\|,\|z-w_2\|\}.$ Let $t'\in(0,1]$ be
the largest number such that $\gamma(t)\in D$ for $t\in(0,t').$ If
$t'=1,$ then an ``integration along $\gamma$'' gives
$$g^\varepsilon_D(z,w_2)-g^\varepsilon_D(z,w_1)\le\pi
c_2\frac{\|w_1-w_2\|}{\min\{\|z-w_1\|,\|z-w_2\|\}}.$$ If $t'<1$,
then, $\gamma(t')\in\partial D.$ Since
$$\lim_{w\to\partial D}g^\varepsilon_D(z,w)=0>
g^\varepsilon_D(z,w_2)>g^\varepsilon_D(z,w_1)$$ and
$g^\varepsilon_D$ is continuous, we way find a $t^\ast\in[0,t')$
with $g^\varepsilon_D(z,\gamma(t^\ast))=g^\varepsilon_D(z,w_2).$
Then, similar as above, we get the same estimates. Letting
$\varepsilon\to 0$ gives the estimate in (iii) with $C=\pi c_2.$

To prove (\ref{8}), we may assume that
$g^\varepsilon_D(z,w_1)<g^\varepsilon_D(z,w_2).$ Let now
$f=f_{z,w}\in\mathcal O(D,\Bbb D)$ be an extremal function for the
Carath\'eodory distance $c_D(z,w)$ (cf. \cite{Jar-Pfl1}, page 16).
We may assume that $f(z)=0.$ For $z\neq w$ set
$h_{z,w}(\zeta)=f(\zeta)/f(w),$ $\zeta\in D.$ Then there are
$c_1,c_3>0$ with
\begin{equation}\label{9}
|h_{z,w}(\zeta)|\le \frac{\tanh c_D(z,\zeta)}{\tanh c_D(z,w)}\leq
c_1\frac{\|\zeta-z\|}{\|w-z\|}\le\frac{c_3}{\|w-z\|},\ z\in K,w\in
D\setminus\{z\}
\end{equation}
(use that  $\Bbb B(z,r)\subset D\subset \Bbb B(z,R)$, $z\in K$,
for certain $r,R$). Set
$$D'=\{\zeta\in D:\zeta+h_{z,w_1}(\zeta)(w_2-w_1)\in D\},\quad
D''=D'\setminus\overline{\Bbb B(z,\varepsilon)}$$ and
$$\hat g(\zeta)=g^\varepsilon_D(z,\zeta+(w_2-w_1)h_{z,w_1}(\zeta)),\quad
\zeta\in D'.$$ It follows by (\ref{7}) and (\ref{9}) that $\Bbb
B(z,\varepsilon)\Subset\Bbb B(z,r/2)\Subset D'$ and $w_1\in D''.$
On the other hand, by (\ref{5}), there is a $c_4>0$ such that
$$g_D(z,\zeta)\ge-c_4\dist(\zeta,\partial D),\quad \zeta\in D
\setminus\Bbb B(z,r/2).$$ This, (\ref{6}) and (\ref{9}) implies
that
$$\min_{\zeta\in\partial D'}
g^\varepsilon_D(z,\zeta)\ge\min_{\zeta\in\partial D'}
g_D(z,\zeta)\ge-c_4\max_{\zeta\in\partial D'}\dist(\zeta,\partial
D)\ge-c_3c_4\frac{\|w_2-w_1\|}{\|w_1-z\|}.$$ Then for
$$v(\zeta)=\hat g(\zeta)-g_D^\varepsilon(z,\zeta)$$
we have that
$$\limsup_{\zeta\to\partial D'}
v(\zeta)\le c_3c_4\frac{\|w_2-w_1\|}{\|w_1-z\|}.$$ On the other
hand, for $\zeta\in\partial\Bbb B(z,\varepsilon),$ it follows by
(\ref{6}) and (\ref{9}) that
$$v(\zeta)\le\log\frac{\max\{\varepsilon,\|\zeta+(w_2-w_1)h_{z,w_1}(\zeta)-z\|\}}{r}
-\log\frac{\varepsilon}{r}$$
$$=\log^+\frac{\|\zeta+(w_2-w_1)h_{z,w_1}(\zeta)-z\|}{\varepsilon}\le
\log\left(1+c_1\frac{\|w_2-w_1\|}{\|w_1-z\|}\right).$$

Since $g_D^\varepsilon(z,\cdot)$ is a maximal plurisubharmonic
function on $D''$ and it is continuous on $\overline{D''}\subset
\overline{D'},$ the domination principle implies that
$$v(\zeta)\le c_2\frac{\|w_2-w_1\|}{\|w_1-z\|},\quad\zeta\in D'',$$
where $c_2=\max\{c_1,c_3c_4\}.$ Applying this for $\zeta=w_1$
gives (\ref{8}).\qed
\smallskip

{\it Proof of Corollary \ref{cor6}.} Recall that there is a $c'>0$
such that $\tilde g_D(z,w)\le c'\|z-w\|,$ $z\in K,w\in D.$
Therefore, we may assume that $w_1,w_2\neq z $ and
$\|z-w_2\|\le\|z-w_1\|.$ Two cases are possible.

Case 1. $|g_D(z,w_1)-g_D(z,w_2)|<1.$ Then
$$|\tilde g_D(z,w_1)-\tilde
g_D(z,w_2)|=\tilde g_D(z,w_2)|\exp(g_D(z,w_1)-g_D(z,w_2))-1|$$
$$<(e-1)c'\|z-w_2\|.|g_D(z,w_1)-g_D(z,w_2)|\le(e-1)c'C\|w_1-w_2\|,$$
where $C$ is the constant from Proposition \ref{5}.

Case 2. $|g_D(z,w_1)-g_D(z,w_2)|\ge 1.$ Then, by Proposition 5,
$$C\|w_1-w_2\|\ge\|z-w_2\|\ge\|z-w_1\|-\|w_1-w_2\|$$
and hence $(C+1)\|w_1-w_2\|\ge\|z-w_1\|.$ It follows that
$$|\tilde g_D(z,w_1)-\tilde g_D(z,w_2)|<\max\{\tilde g_D(z,w_1),\tilde
g_D(z,w_2)\}$$
$$\le c'\|z-w_1\|\le c'(C+1)\| w_1-w_2\|.$$
\qed

{\it Proof of Proposition \ref{pr7}.} By ({\ref{5}), we may find
$c>0$ such that if
$$D_{z,\varepsilon}=\{u\in D:\tilde g_D(z,u)<\varepsilon\},\quad
z,\in K,\;\varepsilon\in(0,1),$$ then
$$\dist(D_{z,\varepsilon},\partial D)\ge c(1-\varepsilon).$$

First, we shall prove (i). In virtue of Corollary \ref{cor6}, it
is enough to find a $c_1>0$ such that
\begin{equation}\label{10} |\tilde
g_D(z_1,w)-\tilde g_D(z_2,w)|\le c_1\|z_1-z_2\|,\quad z_1,z_2\in
K,w\in D.
\end{equation}
We may assume that $K$ is the closure of a smooth domain. Then
there is a $c_2>0$ such that for any $z_1,z_2$ there is a smooth
curve $\gamma$ in $K$ joining $z_1$ and $z_2$ with $l(\gamma)\le
c_2\|z_1-z_2\|.$ Set
$\varepsilon=\varepsilon_{z_1,z_2}=1-\|z_1-z_2\|/c.$ Then, by
``integration along $\gamma$", it suffices to prove (\ref{10}),
fixing $w$ and assuming that $1-\|z_1-z_2\|/c>\sup_{z\in K}\tilde
g_D(z,w).$ Since $w\in D_{z_1,\varepsilon}$ and
$$D\supset\tilde D=\{z+z_2-z_1:z\in D_{z_1,\varepsilon}\},$$ then
$$\tilde g_D(z_1,w)=\varepsilon\tilde g_{D_{z_1,\varepsilon}}(z_1,w)=
\varepsilon\tilde g_{\tilde D}(z_2,w+z_2-z_1)\ge\tilde
g_D(z_2,w+z_2-z_1)$$ (cf. \cite {Zwo}, Lemma 4.2.7, for the first
equality). Hence
$$\tilde g_D(z_2,w)-\tilde g_D(z_1,w)\le\tilde g_D(z_2,w)-\varepsilon
\tilde g_D (z_2,w+z_2-z_1)$$
$$\le C\|z_2-z_1\|+(1-\varepsilon)=(C+1/c)\|z_2-z_1\|,$$ where $C$ is
the constant from Corollary \ref{cor6}. By symmetry,
$$\tilde g_D(z_1,w)-\tilde g_D(z_2,w)\le (C+1/c)\|z_2-z_1\|$$
which implies (\ref{10}).

To prove (ii), it is enough to show that:

$\bullet$ there is a $c_3>0$ such that for any $X,Y\in\Bbb C^n,$

\begin{equation}\label{11}
|A_D(z;X)-A_D(z;Y)|\le c_3\|X-Y\|;
\end{equation}

$\bullet$ if $c_4=\max_{z\in K,\|Z\|=1}A_D(z;Z),$ then
\begin{equation}\label{12}
|A_D(z_1;X)-A_D(z_2;X)|\le c_4\|X\|.\|z-w\|/c
\end{equation}
for any $z_1,z_2\in K$ with $\varepsilon=\varepsilon_{z_1,z_2}>0$
and any $X\in\Bbb C^n.$

Observe that (\ref{11}) follows by choosing $c_3$ such that
$$|\tilde g_D(z,w_1)-\tilde g_D(z,w_2)|\le c_3\|w_1-w_2\|,\quad
z\in K,w_1,w_2\in D$$ and using that hyperconvexity implies
$$A_D(z;X)=\lim_{t\nrightarrow 0}\frac{\tilde
g_D(z,z+tX)}{|t|}.$$

To show (\ref{12}), we may assume that $A_D(z_1;X)\le A_D(z_2;Y).$
Since
$$A_D(z_1;X)=\varepsilon A_{D_{z_1,\varepsilon}}(z_1;X)=
\varepsilon A_{\tilde D}(z_2;X)\ge\varepsilon A_D(z_2;X)$$ (cf.
\cite {Zwo}, Lemma 4.2.7, for the first equality), then
$$0\le A_D(z_2;X)-A_D(z_1;X)\le(1-\varepsilon)A_D(z_2;X)$$
which implies (\ref{12}).\qed

\end{document}